\documentclass{amsart}
\usepackage{amssymb,latexsym}
\usepackage{amsmath}
\usepackage{graphicx}
\usepackage{enumerate}
\usepackage{newlattice}
\usepackage{gensymb}

 \newcommand{\seven}{\SfS 7}
 
\theoremstyle{plain}
\newtheorem{theorem}{Theorem}
 
 \theoremstyle{definition}
 \newtheorem{definition}[theorem]{Definition}
 
\begin{document}
\title[On $\mathcal{C}_1$-diagrams]
{Notes on planar semimodular lattices. IX.\\
$\mathcal{C}_1$-diagrams}  

\author[G.\ Gr\"atzer]{George Gr\"atzer}
\email{gratzer@me.com}
\urladdr{http://server.maths.umanitoba.ca/homepages/gratzer/}
\address{University of Manitoba}
\date{June 15, 2021}

\begin{abstract} 
A planar semimodular lattice $L$ is \emph{slim} 
if $\mathsf M_3$ is not a sublattice of~$L$. 
In a recent paper, G. Cz\'edli introduced a very powerful diagram type
for slim, planar, semimodular lattices, the $\mathcal{C}_1$-diagrams. 
This short note proves the existence of such diagrams.
\end{abstract}



\maketitle    

\subsection*{Background}

\begin{figure}[h!]
\centerline{\includegraphics[scale = 1]{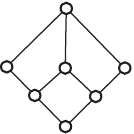}}
\caption{The lattice $\SfS 7$}
\label{F:N7}
\end{figure}

The basic concepts and notation not defined in this note 
are available in Part~I of the book \cite{CFL2}, see \\
{\tt arXiv:2104.06539}\\
it is freely available. 
We will reference it, for instance, as [CFL2, p.~ 4]. 
In~particular, a planar semimodular lattice $L$ is \emph{slim} 
if $\SM{3}$ is not a sublattice of~$L$ and a~grid~$G$ 
is a direct product of two nontrivial chains.
For the lattice $\SfS 7$, see Figure~\ref{F:N7} and  \cite[pages xxi, 34]{CFL2}.
Following my paper  \cite{GKn09} with E. Knapp, a semimodular lattice L is \emph{rectangular} 
if the left and right boundary chains have exactly one doubly-irreducible element each
and these elements are complementary.

In my paper \cite{GLS98a} with H. Lakser, and E.\,T. Schmidt,
we prove that every finite distributive lattice $D$ can be represented
as the congruence lattice of a (planar) semimodulare lattice $L$.
Since $\SM{3}$ sublattices play a crucial role in the construction of $L$, 
it was natural to raise the question
what can be said about congruence lattices of slim, planar, semimodular (SPS) lattices
(see  [CFL2, Problem~24.1], originally raised in my paper~\cite{gG16}). 
The papers in the References list some contributions to this topic.
In particular, my presentation \cite{gG21a} gently reviews
the background of this topic.

\subsection*{$\E C_1$-diagrams}
This research tool played an important role in some recent papers,
see G. Cz\'edli \cite{gC17} and \cite{gCa}, 
G. Cz\'edli and G.~Gr\"atzer~\cite{CG21}, and G.~Gr\"atzer~\cite{gG21a};
for the definition, see G.~Cz\'edli \cite[Definition 5.3] {gC17}, 
G. Cz\'edli~\cite[Definition 2.1]{gCa}, 
and G.~Cz\'edli and G.~Gr\"atzer~\cite[Definition 3.1]{CG21}.

In the diagram of an SPS lattice $K$,
a \emph{normal edge} (\emph{line}) has a slope of $45\degree$ or~$135\degree$.
If it is the first, we call the edge (line) \emph{normal-up}, 
otherwise,  \emph{normal-down}.
Any edge (line) of slope strictly between $45\degree$ and $135\degree$ is \emph{steep}.

A \emph{cover-preserving}~$\SfS 7$ of a lattice $L$ is a sublattice isomorphic to  $\SfS 7$ 
such that the covers in the sublattice are covers in the lattice $L$.

\begin{definition}
A diagram of an SPS lattice $L$ is a
\emph{${\E C}_1$-diagram} if the middle edge of any cover-preserving $\SfS 7$ is steep
and all other edges are normal.
\end{definition}
G. Cz\'edli \cite[Definition 5.11]{gC17} also defines  the much smaller class of $\E C_2$-diagrams.

This note presents a short and direct proof of the existence theorem of ${\E C}_1$-diagrams, 
see G. Cz\'edli \cite[Theorem 5.5]{gC17}, 
utilizing only Theorem~\ref{T:Structure},  the Structure Theorem of Slim Rectangular Lattices.

\begin{theorem}\label{T:well}
Every slim, planar, semimodular lattice $L$ has a ${\E C}_1$-diagram.
\end{theorem}

\begin{figure}[t!]
\centerline{\includegraphics[scale = 1.4]{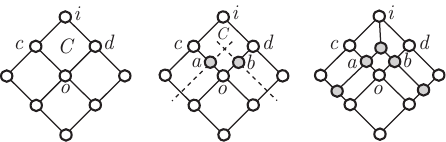}}
\caption{(i) The $4$-cell with $0_C = o$ and  $1_C = i$. 
(ii) Adding the elements $a$ and $b$ for the fork. 
(iii) Adding the fork}
\label{F:notation}
\end{figure}

For an SPS lattice $K$ and $4$-cell $C$ in $K$, 
we denote the \emph{fork extension} of $K$ at $C$ by $K[C]$, 
see G. Cz\'edli and E.\,T.~Schmidt \cite{CS13}
(see also  [CFL2,~Section~4.2]), illustrated by Figure~\ref{F:notation}.

\begin{theorem}[Structure Theorem of Slim Rectangular Lattices]\label{T:Structure}
For every slim rectangular lattice~$K$, there is a grid~$G$ and sequences 
\begin{equation}
G = K_1, K_2, \dots, K_{n-1}, K_n = K\label{E:L}
\end{equation}
of slim rectangular lattices and 
\begin{equation}
C_1 = \set{o_1, c_1,d_1, i_1}, C_2  = \set{o_2, c_2,d_2, i_2} , 
\dots, C_{n-1} = \set{o_{n-1}, c_{n-1},d_{n-1}, i_{n-1}}\label{E:C}
\end{equation}
of $4$-cells in the appropriate lattices such that 
\begin{equation}
G = K_1, K_1[C_1] = K_2, \dots, K_{n-1}[C_{n-1}] = K_n= K \label{E:K}.
\end{equation}
Moreover, the principal ideals $\id c_{n-1}$ and $\id d_{n-1}$ are distributive.
\end{theorem}

\begin{proof}[Proof of Theorem~\ref{T:well} for rectangular lattices]
Let the rectangular lattice $K$ be represented as in \eqref{E:K}.
We prove the Theorem by induction on $n$.
For \text{$n = 1$,}  the statement is trivial.
Let us assume that the statement holds for $n - 1$ and so 
$K_{n-1}$ has ${\E C}_1$-diagrams; we fix one.
By~the induction hypothesis, the $4$-cell $C = C_{n-1}$ with $0_C = o$ and  $1_C = i$
has (at least) two normal edges:~$[o, c]$ and $[o, d]$, see Figure~\ref{F:notation}(i)
and by the last clause of Theorem~\ref{T:Structure},
the principal ideals $\id c$ and $\id d$ are distributive.

Utilizing that $\id c$ is distributive, we place the element~$a$ inside the edge~$[o, c]$
so that the area bounded by the (dotted) normal-up line through $a$
and the normal-up line through $o$ contains no element below $a$;
we place the element~$b$ symmetrically on the other side, 
as in Figure~\ref{F:notation}(ii). 
The two dotted lines meet inside~$C$ 
since the two lower edges of $C$ are normal
and the upper edges are normal or steep.
We~place the third element of the fork at their intersection
and connect it with a steep edge to the element $i$.
We add more elements to the lower left and lower right of $C$ as part of the fork construction,
see Figure~\ref{F:notation}(iii). We can use normal edges for this 
because of the way $a$ and $b$ were placed.
The diagram we obtain is a ${\E C}_1$-diagram of~$K$.
\end{proof}

Now let $K$ be an SPS lattice. 
G. Cz\'edli and E.\,T. Schmidt define in \cite{CS13} a~\emph{corner} element $a$ of~$K$
as a doubly irreducible element on the boundary of $K$ 
such that~$a_*$ is meet-reducible, $a^*$ is join-reducible, 
and $a^*$ has exactly two lower covers.

By G. Cz\'edli and E.\,T. Schmidt \cite{CS13}, 
$K$ is obtained from a slim rectangular lattice~$\hat K$ 
with a fixed ${\E C}_1$-diagram by removing corners.
In a cover-preserving sublattice $\seven$ of $K$, there are only two doubly irreducible elements
but neither is a~corner (since the upper cover of a~corner has
at most two lower covers). Hence, when $\seven$ is a cover-preserving sublattice
(of~$\hat K$ or any other SPS lattice), then this~$\seven$ contains no
corner of K. So the~$\seven$-s remain~$\seven$-s, the steep edges remain the
``legitimately'' steep edges of these remaining~$\seven$-s. All other edges that
are left after removing corners remain of normal slopes. Thus, $K$ is
a ${\E C}_1$-diagram, as required.

\end{document}